\documentclass[12pt]{iopart}
\usepackage{mathptmx}

\newcommand{\Uhsu}{{\mathcal{U}_\hbar(\mathrm{su}_2)}}
\newcommand{\Uhg}{{\mathcal{U}_\hbar(\mathbf{g})}}
\newcommand{\Ug}{{\mathcal{U}\!(\mathbf{g})}}

\newcommand{\UhslC}{{\mathcal{U}_{\hbar}(\mathrm{sl}_2(\mathbf{C}))}}

\newcommand{\slC}{{\mathrm{sl}_2(\mathbf{C})}}
\newcommand{\so}{{\mathrm{so}_4}}
\newcommand{\Uhso}{{\mathcal{U}_{\hbar}(\mathrm{so}_4)}}

\newcommand{\CGs}[6]{
  \setlength{\arraycolsep}{1pt} \Big(
  \begin{array}{cc}
    \scriptstyle #1\! & \scriptstyle #2\vphantom{\scriptstyle #3} \\
    \scriptstyle #4\! & \scriptstyle #5\vphantom{\scriptstyle #6}
  \end{array} \Big|
  \begin{array}{c}
    \vphantom{\scriptstyle #1 #2} \scriptstyle #3 \\
    \vphantom{\scriptstyle #4 #5} \scriptstyle #6
  \end{array}\Big)}

\newcommand{\CGqs}[6]{\CGs{#1}{#2}{#3}{#4}{#5}{#6}_{\!\!\!q}}

\newcommand{\qBinom}[2]{
  \setlength{\arraycolsep}{1pt}\Bigl[
  \begin{array}{c}
    \scriptstyle #1 \\
    \scriptstyle #2
  \end{array} \Big]}

\begin{document}
  
\title{Realization of \textit{q}-deformed spacetime as star product by
  a Drinfeld twist}

\author{Christian Blohmann}

\address{Ludwig-Maximilians-Universit\"at M\"unchen, Sektion Physik\\
Lehrstuhl Prof.\ Wess, Theresienstr.\ 37, D-80333 M\"unchen\\[1em]}

\address{Max-Planck-Institut f\"ur Physik, 
        F\"ohringer Ring 6, D-80805 M\"unchen}

\begin{abstract}
  Covariance ties the noncommutative deformation of a space into a
  quantum space closely to the deformation of the symmetry into a
  quantum symmetry. Quantum deformations of enveloping algebras are
  governed by Drinfeld twists, inner automorphisms which relate the
  deformed to the undeformed coproduct. While Drinfeld twists
  naturally define a covariant star product on the space algebra, this
  product is in general not associative and does not yield a quantum
  space. It is reported that, nevertheless, there are certain Drinfeld
  twists which realize the quantum plane, quantum Euclidean 4-space,
  and quantum Minkowski space.
\end{abstract}

\section{Introduction}

From the beginnings of quantum field theory it had been argued that
the pathological ultraviolet divergences should be remedied by
limiting the precision of position measurements. This is one of the
main motivations to study noncommutative geometries, which imply a
space uncertainty in a natural and fundamental way. From experience we
know that, if spacetime is noncommutative, the noncommutativity can
only be small. This suggests to describe noncommutative spacetime as
perturbative deformation of ordinary, commutative Minkowski space. The
algebraic aspects of a deformation can be separated from the analytic
questions of continuity and convergence by considering formal power
series. In such a framework a noncommutative geometry is a formal
deformation in the sense of Gerstenhaber \cite{Gerstenhaber:1964} of
the function algebra on the space manifold.  Such formal deformations
have appeared naturally in the context of gauge theories on
noncommutative spaces \cite{Seiberg:1999,Madore:2000b}.

Algebraically, physical spacetime is characterized by the Minkowski
algebra of spacetime functions and covariance with respect to the
Lorentz symmetry. The symmetry distinguishes Minkowski space from
other 4-dimensional flat spaces, such as Euclidean 4-space. Covariance
ties the deformation of the symmetry closely to the deformation of the
space. Quantum deformations of the enveloping algebra which describes
this symmetry are known to be governed by Drinfeld twists, inner
automorphisms which relate the deformed to the undeformed coproduct
\cite{Drinfeld:1989,Drinfeld:1989b}. Therefore, one ought to be able
to use these twists in order to deform the space algebra into the
according quantum space as it was suggested in \cite{Grosse:2001}. It
will be shown that for quantum Minkowski space this is indeed
possible.

\section{The problem}

\subsection{Covariant quantum spaces}

Let $\mathbf{g}$ be the Lie algebra of the symmetry group of a space
and $\mathcal{X}$ be the function algebra of this space. The elements
$g \in \mathbf{g}$ of the Lie algebra act on $\mathcal{X}$ as
derivations, $g \triangleright xy = (g\triangleright x)y + x(g
\triangleright y)$ for $x,y \in \mathcal{X}$. A generalized way of
writing this is
\begin{equation}
\label{eq:ModuleAlgebra}
  g \triangleright xy =
  (g_{(1)} \triangleright x)(g_{(2)} \triangleright y)
\end{equation}
for all $g \in \Ug$, where the coproduct is defined as
$g_{\scriptscriptstyle(1)} \otimes g_{\scriptscriptstyle(2)} \equiv
\Delta(g) := g \otimes 1 + 1 \otimes g$ on the generators $g \in
\mathbf{g}$ and extended to a homomorphism on the enveloping algebra.
This covariance condition makes sense for the action of a general Hopf
algebra on an associative algebra $\mathcal{X}$. It ties a deformation
$\mu \rightarrow \mu_\hbar$ of the multiplication map $\mu(x \otimes
y) := xy$ of the space algebra to the deformation $\Delta \rightarrow
\Delta_\hbar$ of the coproduct of the symmetry Hopf algebra if
covariance is to be preserved,
\begin{equation}
\label{eq:BothDeform}
  g \triangleright xy =
    (g_{(1)} \triangleright x)(g_{(2)} \triangleright y) 
    \quad\stackrel{\hbar}{\longrightarrow}\quad
  g \triangleright (x \star y) =
    (g_{(1_\hbar)} \triangleright x)\star (g_{(2_\hbar)}
    \triangleright y) \,,
\end{equation}
where
\begin{equation}
  x \star y := \mu_\hbar(x \otimes y) \quad\mathrm{and}\quad
  g_{(1_\hbar)} \otimes g_{(2_\hbar)} := \Delta_\hbar(g) \,.
\end{equation}
A large class of deformations which are covariant in this sense are
quantum spaces.

\subsection{Star products by Drinfeld twists}
\label{sec:Quest}

In the case of quantum spaces, the deformed coproduct belongs to the
Drinfeld-Jimbo deformation $\Uhg$ of the enveloping Hopf algebra
\cite{Drinfeld:1985,Jimbo:1985}. Drinfeld has observed
\cite{Drinfeld:1989,Drinfeld:1989b} that as $\hbar$-adic algebras
$\Uhg$ and $\Ug[[\hbar]]$ are isomorphic and that the deformed
coproduct $\Delta_\hbar$ of the Hopf algebra $\Uhg \cong
(\Ug[[\hbar]], \Delta_\hbar, \varepsilon_\hbar, S_\hbar)$ is related
to the undeformed coproduct $\Delta$ by an inner automorphism. That
is, there is an invertible element $\mathcal{F} \in (\Ug \otimes
\Ug)[[\hbar]]$ with $\mathcal{F} = 1 \otimes 1 + \mathcal{O}(\hbar)$,
called Drinfeld twist, such that
\begin{equation}
\label{eq:StarProd4}
  \Delta_\hbar(g) = \mathcal{F}\Delta(g)\mathcal{F}^{-1} \,.
\end{equation}
Comparing the covariance condition~(\ref{eq:BothDeform}) of the
deformed multiplication,
\begin{equation}
\label{eq:StarProd3}
  g \triangleright \mu_\hbar(x \otimes y) =
  \mu_\hbar(\Delta_\hbar(g)\triangleright [x \otimes y])
  = \mu_\hbar(\mathcal{F}\Delta(g)\mathcal{F}^{-1}
  \triangleright [x \otimes y])
\end{equation}
with the covariance property~(\ref{eq:ModuleAlgebra}) of the
undeformed product, we see that Eq.~(\ref{eq:StarProd3}) is naturally
satisfied if we define the deformed product by
\begin{equation}
\label{eq:StarProd2}
  \mu_\hbar(x \otimes y)
  := \mu(\mathcal{F}^{-1} \triangleright [x \otimes y])
  \quad \Leftrightarrow \quad
  x \star y :=
  (\mathcal{F}^{-1}_{[1]} \triangleright x)
  (\mathcal{F}^{-1}_{[2]} \triangleright y)\,,
\end{equation}
as it was observed in \cite{Grosse:2001} (suppressing in a Sweedler
like notation the summation of $\mathcal{F} = \sum_i
\mathcal{F}_{\scriptscriptstyle 1i} \otimes
\mathcal{F}_{\scriptscriptstyle 2i} \equiv
\mathcal{F}_{\scriptscriptstyle[1]} \otimes
\mathcal{F}_{\scriptscriptstyle[2]}$). Since the elements of the Lie
algebra $\mathbf{g}$ act on the undeformed space algebra $\mathcal{X}$
as derivations, $\mathcal{F}^{\scriptscriptstyle -1}$ acts as
$\hbar$-adic differential operator on $\mathcal{X} \otimes
\mathcal{X}$. Hence, writing out the $\hbar$-adic sum of
$\mathcal{F}^{\scriptscriptstyle -1} = 1 \otimes 1+ \sum_k \hbar^k
\mathcal{F}_k^{\scriptscriptstyle -1}$ we can define the
bidifferential operators
\begin{equation}
\label{eq:Bdef}
  B_k(x,y) := \mu(\mathcal{F}_k^{-1} \triangleright [x \otimes y])
  = (\mathcal{F}^{-1}_{k[1]} \triangleright x)
  (\mathcal{F}^{-1}_{k[2]} \triangleright y)\,,
\end{equation}
such that the star product~(\ref{eq:StarProd2}) can be written in the
more familiar form \cite{Bayen:1978a}
\begin{equation}
\label{eq:StarProd1}
  x \star y := xy + \hbar B_1(x,y) + \hbar^2 B_2(x,y) + \ldots 
\end{equation}

\subsection{The problem of associativity}

Even though the twist $\mathcal{F}$ yields by Eq.~(\ref{eq:StarProd4})
a coassociative coproduct, Eq.~(\ref{eq:StarProd2}) will in general
not define an associative product. The associativity condition $(x
\star y) \star z = x \star ( y \star z )$ for $\mu_\hbar$ can be
expressed with the Drinfeld coassociator
\begin{equation}
\label{eq:CoassDef}
  \Phi := (\Delta \otimes \mathrm{id})(\mathcal{F}^{-1})\,
  (\mathcal{F}^{-1} \otimes 1) \,(1 \otimes \mathcal{F})
  \,(\mathrm{id} \otimes \Delta)(\mathcal{F}) \,,
\end{equation}
\pagebreak as
\begin{equation}
\label{eq:Quest1}
  (\Phi_{[1]} \triangleright x) (\Phi_{[2]} \triangleright y)
  (\Phi_{[3]} \triangleright z) = xyz \,.
\end{equation}
for all $x,y,z \in \mathcal{X}$.

For a given $\Uhg$-covariant quantum space, is there a Drinfeld twist
$\mathcal{F}$ which yields by Eq.~(\ref{eq:StarProd2}) the associative
product of the quantum space? We will answer this question positively
for three important cases: the quantum plane, quantum Euclidean
4-space, and quantum Minkowski space.

\section{Constructing covariant star products}
\label{sec:StarProducts}

\subsection{The general approach}

To our knowledge, no Drinfeld twist for the Drinfeld-Jimbo quantum
enveloping algebra of a semisimple Lie algebra has ever been computed.
This indicates that it will be rather difficult to answer this
question on an algebraic level. The representations of Drinfeld
twists, however, can be expressed by Clebsch-Gordan coefficients
\cite{Curtright:1991,Blohmann:2002a}. Therefore, we propose the
following approach, which tackles the problem on a representation
theoretic level:

Consider a $\Uhg$-covariant quantum space algebra $\mathcal{X}_h$ and
its undeformed limit, the $\Ug$-covariant space algebra $\mathcal{X}$.

\begin{enumerate}
  
\item Determine the irreducible highest weight representations of all
  possible Drinfeld twists from $\Delta$ to $\Delta_\hbar$.
  
\item Determine the basis $\{ T_{m}^j \}$ of the quantum space
  $\mathcal{X}_\hbar$ which completely reduces $\mathcal{X}_\hbar$
  into (possibly degenerate) irreducible highest weight-$j$
  representations of $\Uhg$.
  
\item Calculate the multiplication map $\mu_\hbar$ of
  $\mathcal{X}_\hbar$ with respect to this basis. The undeformed limit
  $\mu = \lim_{\hbar \rightarrow 0} \mu_\hbar$ yields the commutative
  multiplication map with respect to this basis.
  
\item Check on the level of representations if one of the twists
  realizes the deformed multiplication by Eq.~(\ref{eq:StarProd2}) as
  linear map with respect to this basis.

\end{enumerate}
Since this procedure reduces the algebraic problem to a representation
theoretic one, it works well for cases where the representation theory
is well understood, such as for the quantum spaces of $\Uhsu$,
$\Uhso$, and $\UhslC$.

\subsection{Example: the quantum plane}

The $\hbar$-adic quantum plane generated by $x$ and $y$ with
commutation relations $xy = q yx$, $q := \mathrm{e}^\hbar$, is a
$\Uhsu$-covariant space.  Let us denote by $\rho^{j}$ the structure
map of the spin-$j$ representation of $\Uhsu$. The results of the
proposed approach are:
\begin{enumerate}
\item The irreducible representations of the Drinfeld twists can be
  expressed by the $q$-deformed and undeformed Clebsch-Gordan
  coefficients \cite{Schmuedgen} as
  \begin{equation}
  \label{eq:TwistRep}
    (\rho^{j_1}\otimes \rho^{j_2})(\mathcal{F})^{m_1 m_2}{}_{m_1' m_2'}
    = \sum_{j,m}\eta(j_1,j_2,j)
    \CGqs{j_1}{j_2}{j}{m_1}{m_2}{m}
    \CGs{j_1}{j_2}{j}{m_1'}{m_2'}{m} \,,
  \end{equation}
  where $\eta(j_1,j_2,j) \in \mathbf{C}[[\hbar]]$ is some complex
  formal power series \cite{Curtright:1991,Blohmann:2002a}.
  
\item A basis of the irreducible spin-$j$ $\Uhsu$-subrepresentation of
  the quantum plane is
\begin{equation}
  T^j_m =
  \qBinom{2j}{j+m}_{q^{-2}}^{\frac{1}{2}}
  x^{j-m} y^{j+m} ,
  \textrm{ where $\qBinom{j}{k}_{q}$ is the $q$-binomial coefficient.}
\end{equation}

\item The multiplication map with respect to this basis is
\begin{equation}
  \mu_\hbar(T^{j_1}_{m_1} \otimes T^{j_1}_{m_1}) =
  \CGqs{j_1}{j_2}{j_1 + j_2}{m_1}{m_2}{m_1 + m_2}
  \,\,T^{j_1+j_2}_{m_1+m_2} \,.
\end{equation}
For the undeformed limit $\mu_\hbar \rightarrow \mu$ the
$q$-Clebsch-Gordan coefficient has to be replaced by the undeformed
Clebsch-Gordan coefficient.

\item The twist $\mathcal{F}$ which yields $\mu_\hbar$
  by~(\ref{eq:StarProd2}) can now be read off using the orthogonality
  of the Clebsch-Gordan coefficients to be the one with
  $\eta(j_1,j_2,j)$ = 1 in Eq.~(\ref{eq:TwistRep}).

\end{enumerate}

\subsection{Quantum Minkowski space}

It can be shown \cite{Blohmann:2002a} that there is also a twist
$\mathcal{F}_{\so}$ of $\Uhso$ which realizes quantum Euclidean
4-space and a twist $\mathcal{F}_{\!\scriptscriptstyle\slC}$ of the
quantum Lorentz algebra $\UhslC$ which realizes quantum Minkowski
space by~(\ref{eq:StarProd2}). These twists are composed out of the
twist $\mathcal{F}$ which realizes the quantum plane and the universal
$\mathcal{R}$-matrix of $\Uhsu$ as
\begin{equation}
  \mathcal{F}_{\so} = \mathcal{F}_{13} \mathcal{F}_{24} \,,\qquad
  \mathcal{F}_{\!\slC} = \mathcal{R}^{-1}_{23}
    \mathcal{F}_{13} \mathcal{F}_{24} \,,
\end{equation}
where we use tensor leg notation, $\mathcal{F}_{13} =
\mathcal{F}_{\scriptscriptstyle[1]} \otimes 1 \otimes
\mathcal{F}_{\scriptscriptstyle[2]} \otimes 1$, etc. These expressions
are plausible, considering the fact that $\Uhso$, is the product of
two copies of $\Uhsu$ and that $\UhslC$ is $\Uhso$ twisted by the
$\mathcal{R}$-matrix.

\section{Conclusion}

By definition, Drinfeld twists yield the deformation of an enveloping
algebra into a quantum enveloping algebra.  It was shown that out of
all twists of the Drinfeld-Jimbo algebras $\Uhsu$ and $\Uhso$, and the
quantum Lorentz algebra $\UhslC$ there are certain twists which
realize their fundamental covariant quantum spaces, the quantum plane,
quantum Euclidean 4-space, and quantum Minkowski space, respectively,
as covariant star products on the undeformed, commutative space
algebras. In other words, these particular twists describe the
deformation of space and symmetry completely.

Therefore, it can be expected that these twists also describe
constructions which are solely based on this deformation, such as the
realization of the quantum Minkowski space algebra within the
undeformed Poincar\'{e} algebra or the formal equivalence of deformed
and undeformed gauge theory which was conjectured by Seiberg and
Witten \cite{Seiberg:1999}.

\end{document}